# The characterization problem for one class high order ordinary differential operators with periodic coefficients.

Efendiev R.F.
Instityte Applied Mathematics.
Z.Khalilov, 23, 370148, Baku, Azerbaijan
rakibaz@yahoo.com

Mathematics Subject classification: 34B25, 34L05, 34L25, 47A40, 81U40

## 1.Introduction

The main purpose of the present work is solving the characterization problem which consist of identification of necessary and sufficient conditions on the scattering data ensuring that the reconstructed potential belongs to a particular class. In our case the potential belongs to $Q_+^2$ consisting of functions

$$p_\gamma(x) = \sum_{n=1}^{\infty} p_{\gamma n} \exp(inx) , \quad \sum_{\gamma=0}^{2m-2} \sum_{n=1}^{\infty} n^\gamma |p_{\gamma n}| < \infty. \quad (1.1)$$

which is subclass of the class $Q^2$ of all $2\pi$ periodic complex –valied funtions is considered on the real axis R , belonding to $L_2[0,2\pi]$

The object under consideration is the operator L , given by the differential expression

$$l(y) = (-1)^m y^{(2m)} + \sum_{\gamma=0}^{2m-2} p_\gamma(x) y^{(\gamma)}(x) \qquad (1.2)$$

in the $L_2(-\infty,\infty)$, with potentials $p_\gamma(x) \in Q^2{}_+$.

Note, that some characterizations for the Sturm-Luivvil operator in the class of real-valued potentials belonding to $L_1^1(R)$ ( $L_\alpha^1(R)$ is the class of measurable potentials satisfinding the condition $\int_R dx (1+|x|)^\alpha |p_\gamma(x)| < \infty$ ), has been given by Melin [9] and Marchenko [8].For more detals reference may be found in [5,7,10].

The inverse problem for the potentials form of (1.1) firstly was formulated and solved in [1],where is shown,that the equation $l(y) = \lambda^{2m} y$, has the solution

$$\varphi(x,\lambda\omega_\tau) = e^{i\lambda\omega_\tau x} + \sum_{j=1}^{2m-1} \sum_{\alpha=1}^{\infty} \sum_{n=1}^{\alpha} \frac{V_{n\alpha}^{(j)}}{n + \lambda\omega_\tau(1-\omega_j)} e^{(i\lambda\omega_\tau + i\alpha)x}, \tau = \overline{0,2m-1},$$

$$\omega_j = \exp(ij\pi/m). \qquad (1.3)$$

and Wronscian of the system of solutions,being equal to $(i\lambda)^{m(2m-1)} A$, where

$$A = \begin{vmatrix} 1 & 1 & .. & 1 \\ \omega_1 & \omega_2 & .. & \omega_{2m-1} \\ .. & .. & .. & .. \\ \omega_1^{2m-1} & \omega_2^{2m-1} & .. & \omega_{2m-1}^{2m-1} \end{vmatrix}$$

is not zero at $\lambda \neq 0$.

The limit

$$\varphi_{nj}(x) \equiv \lim_{\lambda \to -\lambda_{nj}} (\lambda + \lambda_{nj})\varphi(x,\lambda), \quad \lambda_{nj} = -\frac{n}{1-\omega_j}, \quad n \in N, j = \overline{1,2m-1}$$

also is a solution of the equation $l(y) = \lambda^{2m} y$, but is lineary dependent with $\varphi(x, \lambda_{nj}\omega_j)$. Thus, the numbers $\widetilde{S}_{nj}, n \in N, j = \overline{1, 2m-1}$, there exist, for which

$$\varphi_{nj}(x) = \widetilde{S}_{nj}\varphi(x, \lambda_{nj}\omega_j), \tag{1.4}$$

Futher as was estalished by Gasymov M.G. [1] if

I. $\sum_{n=1}^{\infty} n |\widetilde{S}_n| < \infty$

II. $4^{m-1} a_m \sum_{n=1}^{\infty} \frac{|\widetilde{S}_n|}{n+1} = p < 1$

where

$$a_m = \max_{\substack{1 \le j \le l \le 2m-1 \\ 1 \le n, r < \infty}} \frac{|(1-\omega_j)(n+r)|}{|r(1-\omega_j) - n(1-\omega_l)\omega_j|} \quad , \quad \widetilde{S}_n = \sum_{j=1}^{2m-1} n^{2m-2} |\widetilde{S}_{nj}| \tag{1.5}$$

then unicually defined functions $p_\gamma(x), \gamma = \overline{0, 2m-2}$ form of (1.1) exist for which the numbers $\{\widetilde{S}_n\}$ are calculated by formula (1.3)- (1.4). Further the full solution of this problem at m=1 was given in [5], where the authors proved the following

Theorem 1: In order to the given sequence of complex numbers be the set of spectral data for the operator $L = \left(\frac{d}{dx}\right)^2 + p_0(x)$, with potential $p_0(x) \in Q_+^2$, it is necessary and sufficient simultaneously satisfing the conditions

1) $\{n\hat{S}_n\}_{n=1}^{\infty} \in l_2$;
2) Infinite determinant

$$D(z) \equiv \left\| \delta_{nk} + \frac{2\hat{S}_k}{n+k} e^{i\frac{n+k}{2}z} \right\|_{n,k=1}^{\infty}$$

exists ($\delta_{mn}$ -Cronekker's symbol), is continuous, is not equal to zero in the closed semi-plane $\overline{C_+} = \{z : \text{Im } z \ge 0\}$ and is analitical inside open semi-plane $C_+ = \{z : \text{Im } z > 0\}$.

In the present work the full inverse problem for the operator (1.2) with potential (1.1) is solved using methods of [2],[5].

Let's formulate the main result of the work.

Definition. Constructed by the help of formulae (1.4) sequence $\{\widetilde{S}_{nj}\}_{n=1, j=1}^{\infty, 2m-1}$, is called as a set of spectral data for the operator (1.2) with potential (1.1)

Thoerem 2. In order to the given sequence of complex numbers $\{\widetilde{S}_{nj}\}_{n=1, j=1}^{\infty, 2m-1}$ to be a spectral data for the operator (1.2) with potential (1.1) it is necessary and sufficient simultaneosly satisfing the following conditions:

1) $\{n\widetilde{S}_n\}_{n=1}^{\infty} \in l_1$ ; \hfill (1.6)

2) ) Infinite determinant

$$D(z) \equiv \det \left\| \delta_{rn} E_{2m-1} + \left\| \frac{i(1-\omega_l)\widetilde{S}_{nj}}{r\omega_l(1-\omega_j) - n(1-\omega_l)} e^{i\frac{n}{1-\omega_j}z} e^{-i\frac{r\omega_l}{1-\omega_l}z} \right\|_{j,l=1}^{2m-1} \right\|_{r,n=1}^{\infty} \tag{1.7}$$

exists ($E_n$-nxn dimensional unit matrix), is continuous, is not equal to zero in the closed semiplane $\overline{C_+} = \{z : \text{Im } z \geq 0\}$ and is analitical inside open semiplane $C_+ = \{z : \text{Im } z > 0\}$.

## 2. On an inverse poblem of the scattering theory in the semiaxis.

On the base of the proof of formulated Theorem 2, as in [5], lies investigation of the equation $l(y) = \lambda^{2m} y$. Taking in it

$$x = it, \lambda = -ik, y(x) = Y(t) \tag{2.1}$$

we get

$$(-1)^m Y^{(2m)}(t) + \sum_{\gamma=0}^{2m-2} Q_\gamma(t) Y^{(\gamma)}(t) = k^{2m} Y(t), \tag{2.2}$$

where

$$Q_\gamma(t) = (-1)^m (-i)^\gamma \sum_{n=1}^{\infty} p_{\gamma n} e^{-nt}, \quad \sum_{\gamma=0}^{2m-2} \sum_{n=1}^{\infty} n^\gamma |p_{\gamma n}| < \infty \tag{2.3}$$

As a result the equation (2.1) is obtained, the potential of which decreases at $t \to \infty$.

<u>Lemma1</u>: The kernel of the operator of transformation (2.2) $K(t,u), u \geq t$, connected with $+\infty$, with potential (2.3) admits representation

$$K(t,u) = \sum_{j=1}^{2m-1} \sum_{n=1}^{\infty} \sum_{\alpha=n}^{\infty} \frac{V_{n\alpha}^{(j)}}{i(1-\omega_j)} e^{\left[-\alpha t + \frac{n}{1-\omega_j}(t-u)\right]}$$

where the series

$$\sum_{j=1}^{2m-1} \sum_{n=1}^{\infty} \frac{1}{n} \sum_{\alpha=n}^{\infty} \alpha^{2m-1}(\alpha - n) \left| V_{n\alpha}^{(j,\tau)} \right|$$

$$\sum_{j=1}^{2m-1} \sum_{\alpha=1}^{\infty} \alpha^{2m-1} \left| V_{\alpha\alpha}^{(j,\tau)} \right|$$

converge.

Proof:

As is shown in [2], the equation (2.2) with potential (2.3) has a solution

$$f(t, k\omega_\tau) = e^{ik\omega_\tau t} + \sum_{j=1}^{2m-1} \sum_{\alpha=1}^{\infty} \sum_{n=1}^{\alpha} \frac{V_{n\alpha}^{(j)}}{in + k\omega_\tau(1-\omega_j)} e^{(ik\omega_\tau - \alpha)t}, \quad \tau = \overline{0, 2m-1}, \tag{2.4}$$

and the numbers $V_{n\alpha}^{(j)}$ are defined by following reccurent formulae:

$$\left[\left(\alpha - \frac{n}{(1-\omega_j)}\right)^{2m} - \left(\frac{n}{(1-\omega_j)}\right)^{2m}\right] V_{n\alpha}^{(j)} = (-1)^{m+1} \sum_{\gamma=0}^{2m-1} \sum_{s=n}^{\alpha-1} \left[i\left(s - \frac{n}{(1-\omega_j)}\right)\right]^\gamma P_{\gamma, s-n} V_{ns}^{(j)} \tag{2.5}$$

at $\alpha = 2, 3, \ldots; n = 1, 2, \ldots, \alpha - 1; j = 1, 2, \ldots, 2m-1$,

$$i^\gamma p_{\gamma\alpha} + \sum_{j=1}^{2m-1} \sum_{n=1}^{\alpha} d_{j\gamma}(n, \alpha) V_{n\alpha}^{(j)} + \sum_{\nu=\gamma+1}^{2m-2} \sum_{j=1}^{2m-1} \sum_{r+s=\alpha} \sum_{n=1}^{s} d_{j\gamma}(n, s, \nu) p_{\nu r} V_{ns}^{(j)} = 0, \tag{2.6}$$

where

$$\frac{1}{n + k(1-\omega_j)} \left[(i\alpha + k)^{2m} - k^{2m} - (i\alpha + k_{nj})^{2m} + k_{nj}^{2m}\right] = \sum_{\gamma=0}^{2m-2} d_{j\gamma}(n, \alpha) k^\gamma; \quad j = \overline{1, 2m-1},$$

$$\frac{(is+k)^\nu - (is+k_{nj})^\nu}{in+k\,(1-\omega_j)} = \sum_{\gamma=0}^{\nu-1} d_{j\gamma}(n,s,\nu)k^\gamma,$$

and the sequence (2.4) admits 2m –times term by term differentiation. Then under the conditions (2.3) we get

$$f(t,k) = e^{ikt} + \int_t^\infty K(t,u)e^{iku}\,du, \qquad (2.7)$$

where

$$K(t,u) = \sum_{j=1}^{2m-1}\sum_{n=1}^{\infty}\sum_{\alpha=n}^{\infty} \frac{V_{n\alpha}^{(j)}}{i(1-\omega_j)} e^{\left[-\alpha t + \frac{n}{1-\omega_j}(t-u)\right]}. \qquad (2.8)$$

Lemma is proved.

Futher, using the methodology of the works [1,2] the equality

$$f_{nj}(t) = S_{nj} f(t, k_{nj}\omega_j), \qquad (2.9)$$

may be obtained, where $f_{nj}(t) = \lim_{k\to k_{nj}}[in + k(1-\omega_j)]f(t,k)$,

$$k_{nj} = -\frac{in}{1-\omega_j},\; j=\overline{1,2m-1}, n\in N.$$

Rewriting (2.9) in the form

$$\sum_{\alpha=n}^{\infty} V_{n\alpha}^{(j)} e^{-\alpha t} e^{\frac{n}{1-\omega_j}t} = S_{nj} e^{\frac{n\omega_j}{1-\omega_j}t} + \sum_{l=1}^{2m-1}\sum_{r=1}^{\infty}\sum_{\alpha=r}^{\infty} \frac{i(1-\omega_j)V_{nr}^{(l)}S_{nj}}{n\omega_j(1-\omega_l)-r(1-\omega_j)} e^{(-\alpha+\frac{n\omega_j}{1-\omega_j})t} \qquad (2.10)$$

multiplying both sides by $\dfrac{1}{i(1-\omega_j)} e^{-\frac{n}{1-\omega_j}u}$ and definig

$$\widetilde{F}(t+u) = \sum_{j=1}^{2m-1}\sum_{n=1}^{\infty} \frac{S_{nj}}{i(1-\omega_j)} e^{\frac{n}{1-\omega_j}(t\omega_j - u)}, t\le u \qquad (2.11)$$

from (2.10) we obtain the Marchenko type equation

$$K(t,u) = \widetilde{F}(t+u) + \int_t^\infty K(t,s)\widetilde{F}(s+u)\,ds. \qquad (2.12)$$

Thus the following Lemma is proved:

<u>Lemma2</u>: If the coefficients $Q_\gamma(t)$ of the equation (2.2) have a form (2.3), then the kernel of the operator of transformation (2.8) satisfyes Marchenko type equation (2.12) at each $t \ge 0$, where the transmission function $\widetilde{F}(t)$ has a form (2.11) and the numbers $S_{nj}$ are defined by (2.9), from which follows that, $S_{nj} = V_{nn}^{(j)}$.

The coefficients $Q_\gamma(t)$ are reconstructed by the kernel of the transformation operator by the help of recurrent formulae (2.5)-(2.6). So, the main equation (2.12) and the form of the transformation operator (2.11) make naturally the formulation of the inverse problem about reconstruction coefficirnts of the equation (2.1) by numbers $S_{nj}$. In this formulation, as in others, in the using of transformatioin operator method, more important moment is the proof of uniqueness solvabitity of the main equation (2.12).

<u>Lemma3</u>. The homogenous equation

$$g(s) - \int_0^\infty \tilde{F}(u+s)g(u)du = 0 \qquad (2.13)$$

with potential $Q_\gamma(t) \in Q_+^2$ has only trivial solution.

Proof: Let $g \in L_2(R^+)$ be a solution of (2.13) and $f$ - be a solution of

$$f(s) + \int_0^s K(t,s)f(t)dt = g(s). \qquad (2.14)$$

Subsituating g into (2.14) and considering (2.12), we get

$$f(s) + \int_0^s K(t,s)f(t)dt + \int_0^\infty [f(u) + \int_0^u K(t,u)f(t)dt]\tilde{F}(u+s)du =$$

$$= f(s) + \int_s^\infty f(t)[\tilde{F}(t+s) + \int_t^\infty K(t,u)\tilde{F}(u+s)du] = 0.$$

From the estimation

$$\left| \tilde{F}(t+s) + \int_0^\infty K(t,u)\tilde{F}(u+s)du \right| \leq Ce^{-s},$$

at $t \geq s$ follows, that $f = 0, g = 0$ and the lemma is proved.

<u>Lemma 4</u>: For each fixed $a, (\operatorname{Im} a \geq 0)$ the homogenous equation

$$g(s) - \int_0^\infty \tilde{F}(u+s-2ia)g(u)du = 0 \qquad (2.15)$$

has only trivial solution in the space $L_2(R^+)$.

Proof:

Replacing x by x+a, where $\operatorname{Im} a \geq 0$, in (1.2) we obtain the equation of same type with potential $Q_\gamma^a(x) = Q_\gamma(x+a)$ satisfing condition(1.1). Note, that the functions $\varphi(x+a, \lambda\omega_j)$ are the solutions of the equations

$$(-1)^m y^{(2m)}(x) + \sum_{\gamma=0}^{2m-2} Q_\gamma^a(x) y^{(\gamma)}(x) = \lambda^{2m} y(x)$$

which at the $x \to \infty$ turns to

$$\varphi(x+a, \lambda\omega_j) = e^{i\lambda\omega_j a} e^{i\lambda\omega_j x} + o(1).$$

So, the functions

$$\varphi^a(x, \lambda\omega_j) = e^{-i\lambda\omega_j a} \varphi(x+a, \lambda\omega_j)$$

also are solutions type of (1.3). Further, define by $S_{nj}(a)$ the spectral data of L with potential $Q_\gamma^a(x)$

$$L \equiv (-1)^m \frac{d^{2m}}{dx^{2m}} + \sum_{\gamma=0}^{2m-2} Q_\gamma^a(x) \frac{d^\gamma}{dx^\gamma}.$$

According to (1.4) we get

$$S_{nj}(a)\, \varphi^a(x, \lambda_{nj}\omega_j) = \lim_{\lambda \to \lambda_{nj}}[n + \lambda(1-\omega_j)]\varphi^a(x,\lambda) = \lim_{\lambda \to \lambda_{nj}}[n + \lambda(1-\omega_j)]\varphi(x+a,\lambda) e^{-i\lambda a} =$$

$$= e^{i\frac{n}{1-\omega_j}a} S_{nj} \varphi(x+a, \lambda_{nj}\omega_j) = e^{i\frac{n}{1-\omega_j}a} e^{-i\frac{n\omega_j}{1-\omega_j}a} S_{nj} \varphi^a(x, \lambda_{nj}\omega_j) = e^{ina} S_{nj} \varphi^a(x, \lambda_{nj}\omega_j)$$

Consequesly

$$S_{nj}(a) = e^{ina} S_{nj}. \qquad (2.16)$$

Now, as above, we get the main equation of the form (2.12) with transmission function

$$\widetilde{F}_a(t+u) = \sum_{j=1}^{2m-1}\sum_{n=1}^{\infty} \frac{S_{nj}(a)}{i(1-\omega_j)} e^{\frac{n}{1-\omega_j}(t\omega_j - u)} = \widetilde{F}(t - ia + u - ia) = \widetilde{F}(t + u - 2ia)$$

From this lemma follows

<u>Theorem 3</u>: The coefficients $Q_\gamma(t)$ of the equation (2.1) satisfying (2.2) are determined unequvocally by numbers $S_{nj}$.

## III. Proof of Theorem 2.

<u>Neccesarity</u>: From the relation (2.9) and form of $f_{nj}(t)$ we obtain

$$S_{nj} = V_{nn}^j.$$

So

$$\sum_{j=1}^{2m-1}\sum_{n=1}^{\infty}n^{2m-1}|S_{nj}| \leq \sum_{j=1}^{2m-1}\sum_{n=1}^{\infty}n^{2m-1}|V_{nn}^{j}| < \infty.$$

i.e. $n^{2m-1}|S_{nj}| \in l_1$. The neccesarity of the condition (1) is proved.

To prove neccesarity of the condition (2) firstly we show, that from the trivial solvability of the main equation (2.12) at t=0 in the class of functions satisfying to inequality $\|g(u)\| \leq Ce^{-\frac{u}{2}}, u \geq 0$, follows trivial solvability in $l_2(1,\infty, R^{2m-1})$ of the infinite system of equations

$$g_{jn} - \sum_{l=1}^{2m-1}\sum_{r=1}^{\infty}\frac{i(1-\omega_j)S_{rl}}{n\omega_j(1-\omega_l)-r(1-\omega_j)}g_{lr} = 0 \tag{3.1}$$

Really, if $\{g_{jn}\} \in l_2, j = \overline{1,2m-1}$, is a solution of this system, then the function

$$g(u) = (g_1(u), g_2(u),...g_{2m-1}(u)) = \sum_{j=1}^{2m-1}\sum_{n=1}^{\infty}S_{nj}g_{jn}e^{-\frac{n}{1-\omega_j}u} \tag{3.2}$$

is defined for all $u \geq 0$, satisfies inequality

$$|g(u)| \leq Ce^{-\frac{u}{1-\omega_j}}; u \geq 0.$$

and is a solution of (2.13), as

$$g(u) - \int_0^{\infty} g(s)\widetilde{F}(u+s)ds = \sum_{j=1}^{2m-1}\sum_{n=1}^{\infty}S_{nj}g_{jn}e^{-\frac{n}{1-\omega_j}u} - \int_0^{\infty}(\sum_{j=1}^{2m-1}\sum_{n=1}^{\infty}S_{nj}g_{jn}e^{-\frac{n}{1-\omega_j}s})(\sum_{l=1}^{2m-1}\sum_{r=1}^{\infty}\frac{S_{rl}}{i(1-\omega_j)}e^{\frac{r}{1-\omega_l}(s\omega_l-u)})ds =$$

$$= \sum_{j=1}^{2m-1}\sum_{n=1}^{\infty}S_{nj}g_{jn}e^{-\frac{n}{1-\omega_j}u} - \sum_{j=1}^{2m-1}\sum_{n=1}^{\infty}\sum_{l=1}^{2m-1}\sum_{r=1}^{\infty}\frac{S_{nj}S_{rl}}{i(1-\omega_l)}g_{jn}e^{-\frac{r}{1-\omega_l}u}\int_0^{\infty}e^{-\frac{n}{1-\omega_j}s}e^{\frac{r\omega_l}{1-\omega_l}s}ds =$$

$$= \sum_{j=1}^{2m-1}\sum_{n=1}^{\infty}S_{nj}g_{jn}e^{-\frac{n}{1-\omega_j}u} - \sum_{j=1}^{2m-1}\sum_{n=1}^{\infty}\sum_{l=1}^{2m-1}\sum_{r=1}^{\infty}\frac{i(1-\omega_j)S_{nj}S_{rl}}{n\omega_j(1-\omega_l)-r(1-\omega_j)}g_{lr}e^{-\frac{n}{1-\omega_j}u} =$$

$$= \sum_{j=1}^{2m-1}\sum_{n=1}^{\infty}S_{nj}e^{-\frac{n}{1-\omega_j}u}[g_{jn} - \sum_{l=1}^{2m-1}\sum_{r=1}^{\infty}\frac{i(1-\omega_j)S_{rl}}{n\omega_j(1-\omega_l)-r(1-\omega_j)}g_{lr}] = 0.$$

So, $g(u) = 0$, therefore, $S_{nj}g_{jn} = 0$ for all, $n \geq 1, j = \overline{1,2m-1}$, and $g_{jn} = 0, j = \overline{1,2m-1}, n \geq 1$ accordingly (3.1). Introduce in the space $l_2(1,\infty; R^{2m-1})$ operator $F(t)$, given by matrix

$$F_{rn}(t) = \|F_{rn}^{jl}\|_{j,l=1}^{2m-1} = \left\|\frac{i(1-\omega_l)S_{nj}}{r\omega_l(1-\omega_j)-n(1-\omega_l)}e^{-\frac{n}{1-\omega_j}t}e^{\frac{r\omega_l}{1-\omega_l}t}\right\|_{j,l=1}^{2m-1} \tag{3.3}$$

and let $\varphi_{2k-1} = \{(\delta_{v1})\delta_{kr}\}_{v,r=1}^{2m-1,\infty}$, $\varphi_{2k} = \{(\delta_{v2})\delta_{kr}\}_{v,r=1}^{2m-1,\infty}$ (($\delta_{ij}$)- column vector), be an orthonormal system in this space. Then, from $n^{2m-1}|S_{nj}| \in l_1$, we get, that

$$\sum_{j,k=1}^{\infty}\left|(F\varphi_j, \varphi_k)_{l_2(1,\infty;R^{2m-1})}\right| < \infty$$

,i.e. $F(t)$ is kernel operator [3]. So there exists the determinant $\Delta(t) = \det(E - F(t))$ of the operator $E - F(t)$ related, as it is not difficult to see, with the determinant $D(z)$ from the condition 2) of the theorem2, by relation $\Delta(-iz) = \det(E - F(-iz)) \equiv D(z)$.

The determinant of the system (3.1) is $D(0)$, and the determinant of analogous system with potential $Q_\gamma^z = Q_\gamma(x+z), \text{Im } z \geq 0$ is

$$D(z) \equiv \det\left\|\delta_{rn} E_{2m-1} - \left\|\frac{i(1-\omega_l)S_{nj}(z)}{r\omega_l(1-\omega_j) - n(1-\omega_l)}\right\|_{j,l=1}^{2m-1}\right\|_{r,n=1}^{\infty} =$$

$$= \det\left\|\delta_{rn} E_{2m-1} - \left\|\frac{i(1-\omega_l)S_{nj}}{r\omega_l(1-\omega_j) - n(1-\omega_l)}e^{inz}\right\|_{j,l=1}^{2m-1}\right\|_{r,n=1}^{\infty} =$$

$$= \det\left\|\delta_{rn} E_{2m-1} - \left\|\frac{i(1-\omega_l)S_{nj}}{r\omega_l(1-\omega_j) - n(1-\omega_l)}e^{i\frac{n}{1-\omega_j}z}e^{-i\frac{r\omega_l}{1-\omega_l}z}\right\|_{j,l=1}^{2m-1}\right\|_{r,n=1}^{\infty}.$$

Therefore to prove the neccesarity of the condition 2) of Theorem 2 it needs to test $\Delta(0) = D(0) \neq 0$. Really, the system (3.1) may be written in $l_2(1,\infty; R^{2m-1})$ as

$$g - F(0)g = 0.$$

As, $F(0)$ is kernel of the operator the Fredholm theory is applicable to it. Accordingly to this theory trivial solvability of it is equaivalent to the fact, that $\Delta(0) = D(0) = \det(E - F(0))$ is not equal to zero [6]. Neccesarity of the condition 2) is proved.

<u>Sufficiently</u>: Let's multiply both sides of the equation (2.12) by $e^{\frac{r\omega_l}{1-\omega_l}u}$ and integrate by $u \in [t, \infty)$. Then we obtain

$$k(t) = F(t)e(t) + k(t)F(t), \qquad (3.4)$$

where $F(t)$ is defined by matrix $\|F_{rn}(t)\|_{r,n=1}^{\infty}$, and

$$e(t) = \|e_{nj}(t)\|_{n,j=1}^{\infty,2m-1} = \left\|e^{\frac{n\omega_j}{1-\omega_j}t}\right\|_{n,j=1}^{\infty,2m-1}; \quad k(t) = \|k_{lr}(t)\|_{l,r=1}^{2m-1,\infty} = \left\|\int_t^\infty K(t,u)e^{\frac{r\omega_l}{1-\omega_l}u}du\right\|_{l,r=1}^{2m-1,\infty};$$

As $F(t)$ is kernel operator for $t \geq 0$ and the condition $\Delta(t) = \det(E - F(t)) \neq 0$
Is satisfied, there exists bounded in $l_2$ inverse operator $R(t) = (E - F(t))^{-1}$.
As $F(t)e(t) \in l_2$, from (3.4) we get

$$k(t) = R(t)F(t)e(t). \tag{3.5}$$

Now, let's take $\langle f, g \rangle = \sum_{n=1}^{\infty} f_n g_n$. Then (2.12) gives

$$K(t,u) = \langle e(t), A(u) \rangle + \langle k(t), A(u) \rangle = \langle e(t), A(u) \rangle + \langle R(t)F(t)e(t), A(u) \rangle = \langle e(t) + R(t)F(t)e(t), A(u) \rangle =$$

$$= \langle R(t)e(t), A(u) \rangle, \tag{3.6}$$

where $A(u)$ is defined as

$$A(u) = \left\| a_{jn}(u) \right\|_{j,n=1}^{2m-1,\infty} = \left\| \frac{S_{nj}}{i(1-\omega_j)} e^{-\frac{n}{1-\omega_j} u} \right\|_{j,n=1}^{2m-1,\infty}$$

Now suppose, the conditions of the theorem are satisfied. Define the function $K(t,u)$ by the equality (3.6) at $0 \le t \le u$ in according to given above considerations. Then by $u \ge t$ we obtain

$$K(t,u) - \int_t^{\infty} K(t,s)F(s+u)ds = \langle R(t)e(t), A(u) \rangle - \int_t^{\infty} \langle R(t)e(t), A(s) \langle e(s), A(u) \rangle \rangle ds =$$

$$= \langle R(t)e(t), A(u) \rangle - \left\langle R(t)e(t), \left\langle \int_t^{\infty} A(s)e(s)ds, A(u) \right\rangle \right\rangle = \langle R(t)e(t), A(u) \rangle -$$

$$- \langle R(t)e(t), \langle A(u), F(t) \rangle \rangle = \langle R(t)e(t), A(u) \rangle - \langle R(t)e(t), A(u)F^*(t) \rangle = \langle R(t)e(t), A(u) - A(u)F^*(t) \rangle =$$

$$= \langle e(t), A(u) \rangle = F(t+u),$$

where "*" denotes the matrix adjoint to $F(t)$ relatively bilinear form $\langle .,. \rangle$. Thus, the following lemma is proved.

**Lemma 5:** At any $t \ge 0$ the kernel $K(t,u)$ of the transformation operator satisfies the main equation

$$K(t,u) = \widetilde{F}(t+u) + \int_t^{\infty} K(t,s)\widetilde{F}(s+u)ds$$

The uniquovacalary solvability of the main equation follows from Lemma 3. By subsituation it is not difficult to calculate, that the solution of the main equation indeed is

$$K(t,u) = \sum_{j=1}^{2m-1} \sum_{n=1}^{\infty} \sum_{\alpha=n}^{\infty} \frac{V_{n\alpha}^{(j)}}{i(1-\omega_j)} e^{\left[-\alpha t + \frac{n}{1-\omega_j}(t-u)\right]},$$

where the numbers $V_{n\alpha}^{(j)}$ are defined from recurrent relations

$$V_{nn}^{(j)} = S_{nj}$$

$$V_{n\alpha+n}^{(j)} = (1-\omega_j)V_{nn}^{(j)} \sum_{l=1}^{2m-1} \sum_{r=1}^{\alpha} \frac{V_{r\alpha}^{(l)}}{r(1-\omega_j) - n(1-\omega_l)\omega_j}.$$

To prove the main statement that the coefficients $Q_\gamma(t)$ have a form (2.2), at first we result the estimations

$$\left| R_{rn}^{jl}(t) \right| \le \delta_{rn}\delta_{jl} + CS_n; \tag{3.8}$$

$$\left|\frac{\partial^{2m-\tau}}{\partial t^{2m-\tau}} R_{rn}^{jl}(t)\right| \le CS_n; \tau = \overline{1,2m-1}, \tag{3.9}$$

for the matrix elements $R_{rn}(t)$ of the operator $R(t)$, there $C = \max\{C_k > 0, k = \overline{1,2m-1}\}$ - ,is constant, $S_n = \sum_{j=1}^{2m-1} n^{2m-1} |S_{nj}|$;

Really from the equality $R(t) = E + R(t)F(t)$ follows

$$\left|R_{rn}^{lj}(t)\right| \le \delta_{rn}\delta_{lj} + \sum_{\tau=1}^{2m-1}\sum_{p=1}^{\infty} \left|R_{rp}^{l\tau}(t)\right|\left|F_{pn}^{\tau j}(t)\right| \le \delta_{rn}\delta_{lj} + 2\sum_{\tau=1}^{2m-1}(\sum_{p=1}^{\infty}\left|R_{rp}^{l\tau}(t)\right|^2)^2(\sum_{p=1}^{\infty}\left|F_{pn}^{\tau j}(t)\right|^2)^2 \le \delta_{rn}\delta_{lj} +$$

$$+ C_1 a_m ((R(t)R^*(t))_{pp} \sum_{p=1}^{\infty} \frac{1}{(n+p)^2}) S_n \le \delta_{rn}\delta_{lj} + C_1 \|R(t)\|_{l_2 \to l_2} S_n.$$

From other hand, as it was noted, the operator-function $R(t) = (E - F(t))^{-1}$ exists and is bounded in $l_2$ (as $F(t)$ is kernel operator by $t \ge 0$ and $\Delta(t) = \det(E - F(t)) \ne 0$) that proves the first inequality of (3.8).

To prove the estimation (3.9) firstly we get the estimation

$$\left|\frac{d}{dt} R_{rn}^{lj}(t)\right| \le \sum_{\tau=1}^{2m-1}\sum_{p,q=1}^{\infty} \left|R_{rp}^{l\tau}(t)\right|\left|\frac{d}{dt} F_{pq}^{\tau k}(t)\right|\left|R_{qn}^{kj}(t)\right| \le \sum_{\tau=1}^{2m-1}\sum_{p,q=1}^{\infty} (\delta_{rp}\delta_{l\tau} + C_2 S_p) S_q (\delta_{qn}\delta_{kj} + C_3 S_n) \le$$

$$\le (1 + C_4 \sum_{p=1}^{\infty} S_p)^2 S_n \le CS_n.$$

Further, considering $\frac{d^l}{dt^l} R(t) = \sum_{n=1}^{l} C_l^n (\frac{d^{l-n}}{dt^{l-n}} R(t))(\frac{d^n}{dt^n} F(t))R(t), l = \overline{1,2m-1}$, by the help of mathematic induction method the inequality (3.9) is proved. In [4] the following relations are proved ( to make correnpondence with our case, let's denote $q_{2m-2-\gamma}(x) = Q_\gamma(x)$).

$$(-1)^m \frac{\partial^{2m}}{\partial x^{2m}} K(x,t) + \sum_{\gamma=0}^{2m-2} q_{2m-2-\gamma}(x) \frac{\partial^\gamma}{\partial x^\gamma} K(x,t) - \frac{\partial^{2m}}{\partial t^{2m}} K(x,t) = 0$$

$$q_0(x) = 2m \frac{d}{dx} K(x,x)$$

$$q_{k+1}(x) = \sum_{\nu=0}^{k} q_\nu(x) \sum_{s=\nu}^{k} C_{2m-3-s}^{k-s} \left\{\frac{\partial^{s-\nu}}{\partial x^{s-\nu}} K(x,t)\Big|_{t=x}\right\}^{(k-s)} +$$

$$+ \sum_{k=0}^{k+2} C_{2m-1-\nu}^{k+2-\nu} \left\{\frac{\partial^\nu}{\partial x^\nu} K(x,t)\Big|_{t=x}\right\}^{(k+2-\nu)} - (-1)^k \frac{\partial^{k+2}}{\partial t^{k+2}} K(x,t)\Big|_{t=x}$$

$$k = 0,1,...2m-3.$$

Now it is not difficult to show,

$$q_0(x) = \sum_{j=1}^{2m-1}\sum_{n=1}^{\infty} \frac{n \cdot S_{nj}}{i(1-\omega_j)} e^{-nt} + \Pi_0(t),$$

where

$$\Pi_0(t) = \sum_{n,p,q,\tau=1}^{\infty} R_{np}^{ej} F_{pq}^{j\tau} e_{q\tau} e_{\tau q} = <R(t)F(t)e(t), A(t)>.$$

is $2i\pi$ periodic function and has bounded (2m-1) order derivatives. It follows from this, that the Fourier coefficients of the function $\Pi_0(-ix), x \in R,$ satisfy $\sum_{n=1}^{\infty} |n^{2m-1}\Pi_n|^2 < \infty$. But, then $\sum_{n=1}^{\infty} n^{2m-2}|\Pi_n| < \infty$. Thus the Fourier coefficients $P_{2m-2,n}$ of the function $Q_{2m-2}(x) = q_0(x)$ satisfy (2.2). Similary for the rest of coefficients $Q_\gamma(x), \gamma = \overline{0, 2m-3}$ it is established, that the Fourier coefficients $p_{\gamma n}$ of the function $Q_\gamma(x) = q_{2m-2-\gamma}(x), \gamma = \overline{0, 2m-3}$ satisfy (2.2). It means, that the Furier coefficients of the function $P_\gamma(x), \gamma = \overline{0, 2m-2}$ satisfy (1.1).

Let, finally, $\{\widetilde{S}_{nj}\}$-be the set of spectral data for the operator $(L - k^{2m}E)$ with constructed coefficients $P_\gamma(x)$. To finish the proof it is enough to show, that $\{S_{nj}\}$ coincides with initial set $\{\widetilde{S}_{nj}\}$. This fact follows from the equality $\widetilde{S}_{nj} = V_{nn}^{(j)} = S_{nj}$. The theorem is proved.